\documentclass[12pt]{article}
\usepackage{amsfonts}
\begin{document}
\def\R{\mathbb{R}}
\def\C{\mathbb{C}}
\def\Z{\mathbb{Z}}
\def\N{\mathbb{N}}
\def\Q{\mathbb{Q}}
\def\D{\mathbb{D}}
\def\Sp{\mathbb{S}}
\def\T{\mathbb{T}}
\def\E{\mathbb{E}}
\def\P{\mathbb{P}}
\def\hb{\hfil \break}
\def\ni{\noindent}
\def\i{\indent}
\def\a{\alpha}
\def\b{\beta}
\def\e{\epsilon}
\def\d{\delta}
\def\De{\Delta}
\def\g{\gamma}
\def\k{\kappa}
\def\qq{\qquad}
\def\L{\Lambda}
\def\E{\cal E}
\def\G{\Gamma}
\def\F{\cal F}
\def\K{\cal K}
\def\A{\cal A}
\def\B{\cal B}
\def\M{\cal M}
\def\P{\cal P}
\def\Om{\Omega}
\def\om{\omega}
\def\s{\sigma}
\def\th{\theta}
\def\Th{\Theta}
\def\z{\zeta}
\def\p{\phi}
\def\m{\mu}
\def\n{\nu}
\def\l{\lambda}
\def\Si{\Sigma}
\def\q{\quad}
\def\qq{\qquad}
\def\half{\frac{1}{2}}
\def\hb{\hfil \break}
\def\half{\frac{1}{2}}
\def\pa{\partial}
\def\r{\rho}
\begin{center}
{\bf GAUSSIAN RANDOM FIELDS: WITH AND WITHOUT COVARIANCES}
\end{center}
\begin{center}
{\bf N. H. BINGHAM and Tasmin L. SYMONS}
\end{center}
\begin{center}
{\it To Mikhailo Iosifovich Yadrenko, on his 90th birthday}
\end{center}
{\bf Abstract} \\
\i We begin with isotropic Gaussian random fields, and show how the Bochner-Godement theorem gives a natural way to describe their covariance structure.  We continue with a study of Mat\'ern processes on Euclidean space, spheres, manifolds and graphs, using Bessel potentials and stochastic partial differential equations (SPDEs).  We then turn from this continuous setting to approximating discrete settings, Gaussian Markov random fields (GMRFs), and the computational advantages they bring in handling large data sets, by exploiting the sparseness properties of the relevant precision (concentration) matrices. \\

\ni {\it Key words} Gaussian random field, covariance function, Bessel potential, stochastic partial differential equation, Mat\'ern process, Gaussian Markov random field, precision matrix, sparseness, numerical linear algebra \\

\ni {\bf 0. Introduction} \\

\i In this survey, we study Gaussian random fields -- Gaussian processes parametrized by some space $\Om$, which will usually be Euclidean space ${\R}^d$ of dimension $d$ or the sphere ${\Sp}^d$ of dimension $d$ as a Riemannian manifold, embedded in ${\R}^{d+1}$; we normalise to radius 1 for convenience.  Gaussian processes are specified by their mean function (which we shall take to be zero for simplicity), and covariance function.  So specifying Gaussian processes is the same as specifying covariance functions, equivalently, positive-definite functions.  All this is in a continuous setting, where we can use our most powerful weapon, calculus. \\
\i Data, however, and the computation needed to handle it, are discrete.  The tension between these two inspired the title of the 2012 paper [SimLR1]: `In order to make spatial statistics computationally feasible, we need to forget about the covariance function' (of course, the authors do not dismiss covariances for theory, only for computational statistics with large data sets).  Regarding Euclidean space and spheres, compare their later comment [SimLP2]: `the great tragedy of spatial statistics is that the Earth turned out not to be flat'. \\
\i We begin in \S 1 with a study of covariance functions in symmetric spaces (this class, aptly named, includes both Euclidean space and spheres).  Our main tool is the powerful and under-utilised {\it Bochner-Godement theorem}.  We continue in \S 2 with {\it Bessel potentials} (again, powerful and under-utilised in this area), stochastic partial differential equations (SPDEs), and Mat\'ern processes, ubiquitous in the field. \\
\i Quoting the above authors again: `In contrast to traditional statistical modelling, practical problems in spatial statistics are, by and large, {\it computational} in nature'.      We turn in \S 3 to Gaussian Markov random fields (GMRFs), which arise [LinRL1] via graphs obtained by triangulation of the space $\Om$.  Dependence occurs only between neighbouring vertices of the graph.  This means that the precision (or concentration) matrix (the inverse of the covariance matrix) is {\it sparse}: most of its elements are zero.  The computational burden is then carried by sparse numerical linear algebra. \\
\i For further background on random fields, see [Mal], [MarP], [Yad]. \\

\ni {\bf 1. Covariances; the Bochner-Godement theorem} \\

\ni {\it 1.1. The Bochner and Bochner-Schoenberg theorems} \\
\i In Euclidean space, the positive definite functions (normalised to be 1 at the origin) are   the characteristic functions, by Bochner's theorem of 1933.  These are also the covariance functions [Kal1, 50]: for a process $\xi$ parametrised by $T$, ${\rho}_{ij} = cov (\xi(\tau_i), \xi(\tau_j))$,
$$
\sum_{ij} a_i a_j \rho_{ij} = \sum a_i a_j cov(\xi(\tau_i), \xi(\tau_j)) = var (\sum_i a_i \xi(\tau_i)) \geq 0,
$$
for all $a_i \in \R, \ \tau_i \in T$.  This result was extended from processes on Euclidean space to isotropic processes on spheres by Bochner and Schoenberg in 1940-42 (see [BinS4] for references).  Here one needs the {\it ultraspherical} (or {\it Gegenbauer}) {\it polynomials} $P_n^{\l}(u)$ of index $\l$ ($u \in [-1,1]$, $\l = \half (d - 1)$); these are classical orthogonal polynomials on the interval $[-1,1]$; see e.g. Szeg\H{o} [Sze, \S 4.7].  The {\it Bochner-Schoenberg theorem} gives the general isotropic covariance on ${\Sp}^d$ to within scale by a {\it convex combination} (or {\it mixture}) {\it of Gegengauer polynomials}:
$$
c \sum_0^{\infty} a_n P_n^{\l}(u), \ \ a_n \geq 0, \ \ \sum_0^{\infty} a_n = 1, \ \ 
u = d({\bf x}, {\bf y}), \ \ {\bf x, y} \in {\Sp}^d, \ \ \l = \half (d - 1)  \eqno(BS)
$$
(here and below, isotropic means that the covariance between the values at two points depends only on their geodesic distance $u = d({\bf x}, {\bf y})$ as here). \\

\ni {\it 1.2.  Symmetric spaces} \\
\i One looks for a general framework to include both these classical results, and this is given by {\it symmetric spaces}.  This field, which stems from Elie Cartan in the 1920s, belongs to differential and Riemannian geometry and Lie theory rather than to probability.  Below we summarise briefly what we need; for background and detail, we refer to the standard works by Helgason [Hel1,2,3,4] and Wolf [Wol1, Ch. 11].  Readers without a geometric background may prefer to think of what are for us the prime examples: spheres, lines and half-lines, and their products.  For the relevant harmonic analysis (extending the Fourier transform to the spherical transform), see \S 1.3, \S 1.4 and \S 1.8 below. \\  
\i A symmetric space is a Riemannian manifold $M$ with curvature tensor invariant under parallel translation.  Equivalently, a {\it geodesic symmetry} exists: this fixes some base point $o$ (the origin in ${\R}^d$, the North Pole in ${\Sp}^d$), and reverses the direction of the geodesics through $o$.  This gives an involutive isomorphism, mapping a point into its `mirror image' on the geodesic through it and $o$.  Then $M$ is a homogeneous space $G/K$: the coset space of $G$, a closed subgroup of the isometry group of $M$ (containing the transvections, [Wol2, 11.1A]) and $K$, the isotropy subgroup of $G$ fixing $o$. \\

\ni {\it 1.3.  The spherical transform and spherical dual} \\
\i To proceed, we need a version of the Fourier transform on symmetric spaces.  This is provided by the {\it spherical transform}, involving spherical measures and spherical functions.  A {\it spherical measure} $m$ is a $K$-bi-invariant multiplicative linear functional on $C_c(K\backslash G/K)$ (continuous functions of compact support on the double-coset space).  The {\it spherical functions} (originally, zonal spherical functions) are the continuous functions $\om: G \to \C$ such that the measure
$$
m_{\om}(f) := \int_G f(g) \om(g^{-1}) d \mu_G(g)
$$
is spherical.  The map
$$
f \mapsto \hat f (\om) := m_{\om}(f) = \int_G f(g) \om(g^{-1}) d \mu_G(g)
$$
is the {\it spherical transform}.  The {\it positive definite} spherical functions $\phi$ on $(G,K)$ are in bijection with the irreducible unitary representations $\pi$ of $G$ with a $K$-fixed unit vector $u$ by
$$
\phi(g) = \langle u, \pi(g) u \rangle.
$$
These form the {\it spherical dual}, $\Lambda$. \\
\i We confine ourselves to the {\it rank-one} case [Hel1, V.6].  These fall into three classes, the compact and non-compact (dual to each other) and the Euclidean (self-dual) [Hel1, V.2].  These are the {\it two-point homogeneous spaces}, the spaces of {\it constant curvature} $\k$: spheres (curvature $\k > 0$), Euclidean ($\k = 0$) and hyperbolic spaces ($\k < 0$) respectively [Hel1, IX.5, X.3], [Wol2]. \\          
\i In the compact rank-one case, the $\pi$ and $\phi$ above, and $\L$, are in bijection with the Cartan-Weyl {\it dominant weights} [Wol1, \S 6.3], and with a subset of $\R$, which we may again write as $\L$, by the Cartan-Helgason theorem [Hel3, V.1.1, 534-538], [Wol1, 11.4B]. \\

\ni {\it 1.4.  The Bochner-Godement theorem} \\
\i That Bochner's theorem may be extended beyond ${\R}^d$ and ${\Sp}^d$ to symmetric spaces goes back to Godement in 1957, giving the Bochner-Godement theorem.  In its modern formulation, this very useful and powerful result is as follows. \\

\ni {\bf Theorem BG (Bochner-Godement theorem)}.  The general positive definite function $\psi$ on a symmetric space $M$ is given (to within scale $c$) by a mixture of positive-definite spherical functions ${\phi}({\l})$ over the spherical dual $\L$ by a probability measure $\mu$:
$$
\psi(m) = c \int_{\L} {\phi}_{\l}(m) \ d \mu(\l) \quad (m \in M): \qquad
\psi = c \int_{\L} {\phi}_{\l} \ d \mu(\l).                                  \eqno(BG)
$$

\i For background and details, see e.g. Wolf [Wol1, Th. 9.3.4], van Dijk [Dij, Prop. 8.3.3] (there as the `Bochner-Godement-Schwartz theorem').  An early use in a probability context (Gaussian processes on compact symmetric spaces) is Askey and Bingham [AskB] in 1976. \\ 
\i Bochner's theorem is the special case $G = {\R}^d$, $K = \{ 0 \}$ (here the spherical functions are the characters, the complex exponentials $e^{ix.} = [t \mapsto e^{ixt}]$).  The Bochner-Schoenberg theorem is the case
$$
{\Sp}^d = SO(d+1)/SO(d) \leftrightarrow (SO(d+1), SO(d)).
$$
Here the spherical functions are the Gegenbauer polynomials above. \\

\ni {\it 1.5.  Products of symmetric spaces; spatio-temporal and geotemporal processes} \\
\i By the Schur product theorem (1911; [Sch], [HorJ, \S 7.5]), products of positive definite functions are positive definite: the class ${\cal{P}}(M)$ of positive definite functions is closed under pointwise products (which correspond to tensor products of representations).  For product spaces, if $f_i \in {\cal{P}}(M_i)$, then $f_1 f_2 \in {\cal{P}}(M_1 \times M_2)$.  The spherical dual of the product is the product of the spherical duals:
$$
M = M_1 \times M_2, \qquad 
(G, K) = (G_1, K_1) \times (G_2, K_2), \qquad 
\L = {\L}_1 \times {\L}_2.
$$
\i The products of interest here are those in which the first factor is spatial (Euclidean or spherical) and the second is time.  We refer to the first case as spatio-temporal and the second as geotemporal.  Both are of great importance, in both theory and applications.  We recall Whittle's advice, that one should think of a spatial process as the equilibnrium process of a spatio-temporal one [Whi3, Preface].  \\

\ni {\it 1.6. Gelfand pairs} \\
\i For $G$ a locally compact group (Lie group will suffice here) with $K$ a compact subgroup, $(G,K)$ is called a {\it Gelfand pair} if the convolution algebra of $K$-bi-invariant continuous measures on $G$ with compact support is commutative.  For symmetric spaces $M = (G,K)$, $(G,K)$ is a Gelfand pair ([Hel1, X, Th. 2.9, 4.1], [Hel3, IV, Th. 3,1]).  We use the language of symmetric spaces rather than Gelfand pairs, but Gelfand pairs (which can be discrete or continuous) are very interesting; see van Dijk [Dij] for a monograph treatment.  Their relevance here is shown by there being a form of the Bochner-Godement theorem for them ([Dij, Prop. 8.3.3] as noted above). \\

\ni {\it 1.7.  Sphere cross line; the Berg-Porcu theorem} \\
\i From {\it 1.5} above, one can produce geotemporal covariances by taking an isotropic covariance on the sphere (from the Bochner-Schoenberg theorem), a stationary covariance on the line (by Bochner's theorem) and taking their product.  This, though a valid covariance on sphere cross line, is {\it separable} (one might as well deal with the two separately).  The interesting and useful covariances here are {\it inseparable}.  The question was raised by Mijatovi\'c and the authors [BinMS, \S 4.4] in 2016 of finding the general (isotropic, stationary) covariance on sphere cross line, and answered by Berg and Porcu [BerP] in 2017. \\
\i We state the Berg-Porcu theorem below, and illustrate the power and usefulness of the Bochner-Godement theorem by using it to give a very short and simple proof (indeed, two proofs). \\

\ni {\bf Theorem BP (Berg-Porcu Theorem)}.  The class of isotropic stationary covariances on sphere cross line coincides (to within scale) with the class of mixtures of products of Gegenbauer polynomials $P_n^{\l}(u)$ and characteristic functions ${\phi}_n(t)$ on the line: with $u$, $\l$ as before, they are given by
$$
c \sum_0^{\infty} a_n P_n^{\l}(u) {\phi}_n(t), \quad 
c > 0, \quad a_n \geq 0, \quad \sum_0^{\infty} a_n = 1.  \eqno(BP)
$$

\ni {\it Proof}.  In the notation above, with the first factor the sphere ${\Sp}^d$ the spherical dual ${\L}_1$ is the set of Gegenbauer polynomials $P_n^{\l}$; with the second factor the line, ${\L}_2$ is the set of characters, which can also be identified with the line:
$$
t \leftrightarrow e^{it.} = (x \mapsto e^{ixt}).
$$
We can now proceed using the language of either measure theory or probability theory. \\
1. {\it Measure theory}.  We use disintegration of measures in $(BG)$ (Fubini's theorem extended beyond product measures: see e.g. [Kal1, Th. 6.4], Bogachev [Bog, \S 10.6]).  Here, ${\L}_1 = {\N}_0 := {\N} \cup \{ 0 \}$, ${\L}_2 = \R$.  Integrate the probability measure $\mu$ on $\L$ over the second (time) variable for fixed $n$, then sum over $n$.  This disintegrates $\mu$ into a sequence of probability measures ${\mu}_n$ on the line, and a probability measure $(a_n)$ on ${\N}_0$.  Now integrating $e^{itx}$ over ${\mu}_n(dx)$ gives its characteristic function ${\phi}_n(t)$, the second factor in $(BP)$.  The remaining integration is a sum over $n$ of the product of this and the first factor $P_n^{\l}(u)$ with weight $a_n$, giving the result. \\
2.  {\it Probability theory}.  Take $\l = ({\l}_1, {\l}_2)$ in $(BG)$ as a random variable with law $\mu$.  Condition on its first coordinate, and use the Conditional Mean Formula [Wil2, 390] (a special case of the tower property [Wil1, 9.7i] or chain rule [Kal1, 105]). 
\hfil $\square$ \break

\ni {\it Products of spheres}.  There is a similar result due to Guella, Menegatto and Peron in 2016 [GueMP].  Its proof is immediate from the Bochner-Godement theorem.  See [BinMS] for details, [BinS3,4] for background.\\

\ni {\it 1.8. Transforms: commutative and non-commutative settings} \\
\i The prototypical transform for us is the Fourier transform.  There are two major extensions to commutative settings: Pontryagin duality on locally compact groups (1934), and the Gelfand transform on commutative Banach algebras in the 1940s [Rud, Ch. 11].  The spherical transform above derives from Gelfand in 1950 [Gel] (see its review by Godement [God]).  The commutativity in a Gelfand pair captures what is needed to extend to a genuinely non-commutative setting. \\     
    
\ni {\bf 2.  Bessel potentials, SPDEs and Mat\'ern processes} \\

\ni {\it 2.1. Riesz potentials} \\
\i The two most lasting contributions of Marcel Riesz (1886-1969, younger brother of Frederick Riez) are Riesz means (typical means) in summability theory, and his work of 1937-38 on fractional potentials [Rie].  As potential theory is intimately linked to the Laplacian, it is natural that one can treat Riesz potentials in terms of fractional Laplacians (an important type of pseudodifferential operators, [Tay]).  See the now classic treatment in Stein [Ste, V].  Now the Laplacian $\De$ has non-positive eigenvalues (as $D^2 e^{i \omega .} = -{\omega}^2 e^{i \omega .}$), so it is preferable to deal with $- \De$, which has non-negative eigenvalues, before taking fractional powers.  The operators 
$$
I_{\a} := (-\De)^{- \a/2}
$$ 
are called {\it Riesz potentials}, and in $d$ dimensions ${\R}^d$ are well-behaved for $\a \in (0,d)$.  Then for suitable functions $f$,
\begin{eqnarray*} 
(-\De)^{- \a /2} (f) 
&=& 
\frac{\Gamma( \frac{d}{2} - \frac{\a}{2})}{{\pi}^{d/2} 2^{\a} \Gamma (\a/2)} 
\int_{{\R}^d} |x - y|^{-d + \a} f(y) dy \quad (0<\a < d) \\
&=& \int_{{\R}^d} k(x,y) f(y) dy, 
\end{eqnarray*}
say, where $k(.,.)$ is called the {\it kernel} of the operator $I_{-\a/2}$ [Ste, V.1]. \\

\ni {\it 2.2.  Bessel potentials} \\
\i The restriction $\a < d$ is often inconvenient [Ste V.3].  We can avoid it by working instead with $(I - \De)$, which has {\it positive} eigenvalues.  For suitable $f$, the Fourier transforms of $-\De f$ and $(I - \De) f$ are $4 {\pi}^2 |x|^2 \hat f (x)$ and $(1 + 4 {\pi}^2 |x|^2) \hat f(x)$. The corresponding fractional powers
$$
{\cal{J}}_{\a} := (I - \De)^{-\a/2}
$$
are called {\it Bessel potentials} (the calligraphic ${\cal{J}}_{\a}$ is to avoid confusion with $J_{\a}$, the Bessel function of the first kind).  The {\it kernel} of ${\cal{J}}_{\a}$ is the function $G_{\a}$ with Fourier transform 
$$
(1 + 4 {\pi}^2 |x|^2)^{-\a/2}.          
$$
The name Bessel potential derives from the fact that the kernel $G_{\a}$ involves a Bessel function of the third kind with imaginary argument, more briefly a {\it Macdonald function} [Wat, 3.7(6)], $K_{\nu}$ say.  With all three kinds of Bessel function, $J_{\nu}, I_{\nu}, K_{\nu}$, the standard notation for the order is $\nu$ [Wat]. \\
\i We now have [Den], [AroS]
$$
G_{\a}(x) :=
\frac{K_{(d-\a)/2}(|x|) \ |x|^{(\a - d)/2}}{2^{(d + \a - 2)/2} {\pi}^{d/2} \ \Gamma(\a/2)}.
$$
The functions in the two last displays are essentially Fourier transforms of each other; see below.\\
\i The dimension $d$ is fixed; the parameter $\a$ governs {\it smoothness}, as we shall see.  It is convenient, following Simpson. Lindgren and Rue [SimLR1], to use 
$$
\eta := (\a - d)/2                                                                   \eqno(eta)
$$
instead as smoothness parameter.  We can then write the kernel $G_{\a}$ as
$$
G_{\a}(x) = |x|^{\eta} K_{\eta}(|x|)/ c_{\eta} , \quad 
c_{\eta} := 2^{(d + \a - 2)/2} {\pi}^{d/2} \ \Gamma(\a/2).
$$
Note that $\eta$ may be negative (not a problem, as $K_{-\eta} = K_{\eta}$; see below). \\
\ni {\it Macdonald functions}. \\
\i As Watson remarks, the importance of the Macdonald function is as a Bessel function with exponential decay at infinity [Wat, 3.7 and 7.23].  The integral representation
$$
K_{\nu}(x) = \frac{1}{2} \int_0^{\infty} u^{\nu} \ \exp  \ \{ x (u + 1/u) \} \ du/u
$$
is often useful [Erd, II, 7.12(23)], [J{\o}r]. \\
\i The behaviour of $K_{\nu}$ near the origin is important, as we shall see.  We note here that $K_0$ has a logarithmic singularity at the origin [Wat, (3.71)(14)]: $K_0(x) \sim \log (1/x)$ ($x \to 0+$).  For $\nu \neq 0$ [Wat, 3.7(2),(6)]:
$$
K_{\nu}(z) :=
\half \pi \ \frac{I_{-{\nu}}(z) - I_{\nu}(z)}{\sin \pi \nu}, \qquad
I_{\nu}(z) = \sum_0^{\infty} \frac{(\half z)^{\nu + 2m}}{m! \Gamma(\nu + m + 1)}.
$$
(using the obvious `L'Hospital convention' for $\nu = 0$), so
$$
K_{\n} = K_{-\nu}.
$$
For small $z$ $K_{\nu}(z) \sim c_{\nu} z^{\nu}$ for constant $c_{\nu} \neq 0$.    \\ 

\ni {\it 2.3. Macdonald-Student pairs} \\
\i Stein gives the following useful integral formula for $G_{\a}$ [Ste, V.3(26)]:
$$
G_{\a}(x) = \frac{1}{(4 \pi)^{\a/2} \Gamma(\a/2)} 
\int_0^{\infty} e^{-\pi |x|^2/u} e^{-u/4\pi} u^{(-d + \a)/2} du/u     \eqno(G_{\a})  
$$
(as he points out, this is just a rephrasing of the definition of the Gamma function).  Thus $G_{\a}$ is positive, even, integrable, and decreasing on $(0,\infty)$.  So it is (to within a scale factor) a probability density, let us call it the {\it Macdonald density}.  So the Macdonald and Student-$t$ densities form a {\it Fourier pair} (`self-reciprocal pair' in Feller's terminology [Fel, 503]): each is the other's characteristic function to within scale (as one is always integrable [Fel, XV.3 Th. 3]).  For this {\it Macdonald-Student Fourier pair} ({\it Macdonald-Student pair} for brevity), see e.g. Guttorp and Gneiting [GutG].  They give a historical study of of the family below, studied in 1960 by the Swedish statistician Bertil Mat\'ern (1921-2010) [Mat], and now named after him. \\
 
\ni {\it 2.4. The Mat\'ern family of correlation functions; Mat\'ern processes} \\
\i Guttorp and Gneiting [GutG, (1)] give the correlation function in the {\it Mat\'ern family} between values of a spatial random function at locations $s \in {\R}^d$ apart as 
$$
{\rho}_{\nu}(s) := 
\frac{2^{1-\nu}}{\Gamma(\nu)} (\k |s|)^{\nu} K_{\nu}(\k |s|)
\propto \int_{{\R}^d} \frac{\exp \{ i s^T x \} }{ ({\k}^2 + |x|^2)^{(2 \nu + d)/2}} \ dx
\eqno(Mat) 
$$
($\k > 0$ is the {\it scale parameter}, $\nu > 0$ is their notation for the {\it smoothness parameter}; the variance is scaled so that ${\rho}_{\nu}(0) = 1$). \\
\i Bochner's theorem of 1933 shows that characteristic functions of probability distributions are exactly the continuous positive definite functions taking the value 1 at the origin.  Khintchine's theorem of 1934 shows that for {\it stationary} processes, the correlation function is the Fourier transform of a probability distribution, the {\it spectral measure}.  More is true: by the {\it Cram\'er representation} of 1942, the process itself is the Fourier transform of a stochastic process with orthogonal increments (for details and references, see e.g. Doob [Doo, X, XI], Cram\'er and Leadbetter [CraL, \S 7.5]).  Thus the Mat\'ern family of correlations are those of stationary processes, the {\it Mat\'ern processes}, with spectral densities of Student type.  Note that here the correlation function does not depend on the dimension $d$, while the spectral density does.   \\
\i Mat\'ern processes have long been widely used, and have become `the workhorses of spatial statistics'.  Their theoretical study is thus justified, not only by its intrinsic interest, but also by the decades of practical experience of practitioners.  This is well summarised by M. L. Stein, whose advice in his book on spatial data and kriging is `use the  Mat\'ern model' [Stei, 14].  He `showed that the behaviour of the covariance function near the origin has fundamental implications on predictive distributions, particularly predictive uncertainty.  The key feature of the Mat\'ern is the inclusion of a smoothness parameter that directly controls correlation at small distances', and `The smoothness parameter is aptly named as it implies levels of mean-square differentiability of the random process, with large $\nu$ yielding very smooth processes that are many times differentiable, and small $\nu$ yielding rough processes'.  There is also a link with Hausdorff dimension [GenK, \S 3]. \\

\ni {\it 2.5. Bessel potentials and the Mat\'ern family} \\
\i The kernel $G_{\a}$ of ${\cal J}_{\a}$ involves a Macdonald function $K_{\eta}$.  Crucial for us is the stochastic partial differential equation (SPDE) 
$$
({\k}^2 - \De)^{\a/2} \ X(s) = \sigma \ W(s) \qquad (s \in {\R}^d),              \eqno(SPDE)
$$
where $W$ is white noise.  We note here that white noise involves generalised functions (Schwartz distributions) rather than ordinary ones.  For background on white noise analysis, we refer to the standard work by Hida et al. [HidKPS], and for background on generalized random fields to Gelfand and Vilenkin [GelV]. \\
\i Following Simpson, Lindgren and Rue [SimLR1], we write the Mat\'ern covariance as
$$
c_{\eta}(s_i, s_j) = C_{\eta}(\Vert s_i - s_j \Vert) 
:= \frac{{\s}^2}{\Gamma(\eta + d/2) (4 \pi)^{d/2} {\k}^{2 \eta} 2^{\eta - 1}}
 \ (\k \Vert s_i - s_j \Vert)^{\eta} K_{\eta}(\k \Vert s_i - s_j \Vert^2)
$$
(this differs slightly from the notation of [GutG]: $\k$ is a scale parameter as there; $\sigma$ appears as this is a covariance, not a correlation; we use $\eta$ for the smoothness parameter here rather than the traditional $\nu$ to show the link with the dimension $d$ and the Bessel parameter $\a$ in (eta) above). \\ 
\i It is clear from the effect of the Bessel potential being to introduce a multiplier into the Fourier transform that a convolution is involved, and clear from the stochastic context that this convolution will involve a stochastic process.  So it is to be expected that the deterministic factor is the kernel $c_{\eta}$ above and the stochastic factor is white noise.  \\   
      
\ni {\it 2.6. Convolution (kernel) representations; Whittle's theorem} \\
 \i The link between Mat\'ern processes and SPDEs goes back to Whittle [Whi1] [Whi2] in 1954 and 1963.  This was the starting-point of the 2011 study [LinRL], and its 2012 sequel [SimLR1].  So we shall call the result below Whittle's theorem; we follow [SimLR1] in the proof.  The result is as in Simpson's contribution to the discussion of [LinRL, 65-66]. \\

\ni {\bf Theorem W (Whittle's theorem)}.  The stationary solutions to $(SPDE)$ above are the convolutions
$$
X(s) = \int_{{\R}^d} c_{\eta}(s,t) \ dW(t) \qquad (s \in {\R}^d)
$$
with $c_{\eta}$ the Mat\'ern covariance above and $W$ white noise.  That is,
$$
\int_{{\R}^d} ds \ \psi(s)X(s) = 
\int_{{\R}^d} ds \ \psi(s) \int_{{\R}^d} c_{\eta}(s,t) \ dW(t)
$$
for all test functions $\psi$. \\

\ni {\it Proof}.  Following Walsh [Wal], we define a solution of $(SPDE)$ to be any random field $X(s)$ satisfying
$$
\int X(s) ({\k}^2 - \De)^{\a/2} \phi(s) ds = \s \int \phi(s) dW(s)
$$
for every test function $\phi$.  Choose $\phi$ to solve
$$
({\k}^2 - \De)^{\a/2} \phi(s) ds = \psi(s):
$$
then  
$$
\phi(s) = ({\k}^2 - \De)^{-\a/2} \psi(s)
$$
is smooth [Ste].  Substituting in the above gives
$$
\int X(s) \psi(s) ds = \s \int ({\k}^2 - \De)^{-\a/2} \psi(s) dW(s).
$$
But as the kernel of the Bessel potential is the Mat\'ern covariance $c_{\eta}$,
$$
({\k}^2 - \De)^{-\a/2} \psi(s) = \int c_{\eta}(s,t) \psi(t) dt.
$$
Substituting and using Fubini's theorem gives the result. \hfil $\square$ \break
   
\i As above, Whittle's result [Whi1], [Whi2] was the starting-point for [LinRL].  The authors proceeded, as in [GutG], to use the above results on Fourier transforms to identify the finite-dimensional distributions of the two sides of $(SPDE)$ by showing that they integrate the same way over each finite set of test functions.  \\
{\it Note}. 1.  One important advantage of the approach above is that the parameter $\a > 0$ is not subject to the restriction to half-integer values in [LinRL], [SimLR1].  Such situations are typical of this area, where the relevant analysis and/or special-function theory involves a continuous parameter, which has geometric significance as a dimension for integer (or half-integer) values.  See e.g. [Bin1,2,3], [BinS1,2]. \\
2.  The representation of the process $X$ in Whittle's Theorem as a convolution of the function (kernel) $C_{\eta}$ and the process $W$ is a prototype of a wider class, the {\it convolution processes} (process convolutions, kernel convolutions).  These were introduced by Higdon [Hig1,2] in 1998, in a climate-modelling setting, and have been widely used since; see e.g. Rodrigues and Diggle [RodD]. \\
3. Kernel convolutions have the great advantage that one can generalise by replacing the Mat\'ern kernel $c_{\eta}$ by any $L_2$-kernel $k$:
$$
X(s) = \int_{{\R}^d} k(s,t) \ dW(t)
$$
still gives a Gaussian random field.  This gives one great flexibility when modelling.  One can also replace white noise by any independently scattered random measure (measures of disjoint sets are independent: [Kal2]), and hence one can construct non-Gaussian random fields [AbeP]. \\              
4. For more on smoothness of paths (regularity) and SPDEs, see [LanS], [BinS1], [BroKLO].   \\

\ni {\it 2.7. Matern processes on spheres, manifolds and graphs} \\
%\i [BorTMD]; functional calculus \\
\i Much of the theory above, motivated by the needs of spatial statistics, has involved Fourier analysis in Euclidean space, which is self-dual under the Fourier transform.  While this is the right framework for a locality or region, it is not on the larger scale (`maps and globes').  We live on Planet Earth, which is (approximately) a sphere, which we may scale to be the unit sphere, $\Sp := {\Sp}^2$ as in \S 1.   \\
\i As on any compact Riemannian manifold $M$, there exists a complete orthonormal basis $\{ {\phi}_n \}$ of $L_2(\Sp)$ of eigenfunctions of $-\De = -{\De}_{\Sp}$, with eigenvalues ${\l}_j$ with
$$
0 = {\l}_0 < {\l}_1 \leq \cdots \leq {\l}_j \cdots \uparrow + \infty
$$
One can expand a suitable function on $M$ in an eigenexpansion, the Sturm-Liouville expansion (see e.g. Chavel [Cha1, VI.1], where this is done for the heat kernel):
$$
- \De f = \sum_0^{\infty} {\l}_n \langle f, {\phi}_n \rangle {\phi}_n.
$$
For suitable $\Phi$, one can extend this to
$$
\Phi(- \De) f = \sum_0^{\infty} \Phi({\l}_n) \langle f, {\phi}_n \rangle {\phi}_n,
$$
an example of the {\it functional calculus} (or {\it symbolic calculus}, of F. Riesz and N. Dunford); see e.g. [Rud, 290-292]. \\ 
\i  To exploit the geometry, one expands in a basis of interest; for  the sphere here, one uses the {\it spherical harmonics} ([AndAR Ch.9], [SteW, IV.2]) as in [BinS1 \S 3]. \\
\i Using $\cal{F}$ for the Fourier transform, the last display in the proof above may be summarised as [LinRL, B.3.2]
$$
\{ {\cal{F}} ((\k^2 - \De)^{\a/2}) \phi \}({\bf k})
= (\k^2 + \Vert {\bf k} \Vert^2)^{\a/2} ({\cal{F}} \phi) ({\bf k}) \qquad ({\bf k} \in {\R}^d).
$$
The analogous statement for the sphere is [LinRL, B.3.2]
$$
\{ {\cal{F}} ((\k^2 - \De)^{\a/2}) \phi \}(k)
= (\k^2 + {\l}_k^2)^{\a/2} ({\cal{F}} \phi) (k) \qquad (k = 0,1,2,\cdots).
$$
For such {\it fractional calculus on spheres}, see [BinS1, \S 4 Remark 2] and the references there (to Askey and Wainger, Bavinck etc.)  \\
\i For Mat\'ern Gaussian processes on general Riemannian manifolds, see Borovitskiy et al. [BorTMD].  For graphs, one proceeds similarly using the {\it graph Laplacian} (see e.g. [Chu]).  For Mat\'ern processes on graphs, see Borovitskiy et al. [BorATMDD, \S 3]. \\

\ni {\it 2.8. Extreme values} \\
\i The area of spatial (or spatio-temporal, or geotemporal) extremes is vast and topical.  For reasons of space, we can only refer here to [AsaDE], [CooNN], [HaaRo], [ShaW] and the first author's survey [BinO] with Ostaszewski. \\ 

\ni {\bf 3. Gaussian Markov random fields; sparse numerical linear algebra} \\
% [LinRL1] FEM \\

\ni {\it 3.1. Precision matrices and sparseness} \\
\i If the space $\Om$ is discretised by triangulation, the vertices and edges of the triangles form a {\it graph}, $G = (V, E)$, $V = \{v_1, \cdots, v_n \}$ say.  We may thus use the powerful {\it graphical models} of statistics; see e.g. Lauritzen [Lau].  \\
\i The values $x_i$ of our Gaussian process $X$ at the vertices $v_i$ are multivariate Gaussian, with covariance matrix $\Si = ({\s}_{ij})$, say, and precision (concentration, inverse covariance) matrix $K := {\Si}^{-1} = (k_{ij})$ (`K for Konzentration').  Two Gaussian variables $x_i$ and $x_j$ are independent if and only if their covariance ${\s}_{ij} = 0$, i.e. they are uncorrelated (for the simple proof, see e.g. [Kal1, Lemma 13.1]). \\
\i That $x_i$, $x_j$ are {\it conditionally independent given all the others} if and only if $k_{ij} = 0$ is much more recent.  Now taking $x_1$ as the random $2$-vector and $x_2$ as the $(p - 2)$-vector we condition on, the Gaussian Regression Formula (GRF; [Lau, Prop. C5], [BinK, 5.2]) shows us that the conditional density $f_{1|2}$ of $x_1 | x_2$ has (conditional) covariance matrix $K_{11}^{-1}$.  So one has conditional independence if and only if $K_{11}^{-1}$ is diagonal, that is (the matrices being $2 \times 2$), $K_{11}$ is diagonal, that is, $k_{12} = 0$. \\
\i   This important result derives from Dempster [Dem2, p.159], [Dem1], so as it needs a name let us call it {\it Dempster's theorem}.  It owes its great importance in modern statistics to the role it plays in {\it graphical models} [Lau, Prop. 5.2].  {\it Sparseness} properties of concentration matrices are highly revealing about structure, as in [HasTW], as well as being numerically very convenient (we note in passing the great importance of numerical linear algebra, here and elsewhere; see e.g. [GolV]).  Conditional independence statements of this sort are important in Markov properties on undirected graphs (Hammersley-Clifford theorem; [Lau, 3.2.1]), multivariate normal models [Lau, 5.1.3], covariance selection models ([Lau, 5.2], following [Dem2]), etc.  They come into their own with Gaussian Markov random fields; see [RueH].\\
\i Building on [RueH], [LinRL] proceed by triangulating the space $\Om$.  For this one can use Delaunay triangulation, giving the dual of the Voronoi diagram (see e.g. [PrepS, Ch. 5,6]).  Here one maximises the minimum angle in a triangle, avoiding `long thin triangles' which are numerically (and geometrically) awkward. \\
\i Although the flexibility given by the smoothness parameter $\nu$ in the Mat\'ern model to           reflect smoothness of the process being modelled is valuable, `$\nu$ usually fixed, since it is poorly identified in typical applications.  A more natural interpretation of the scale parameter $\k$ is as a {\it range parameter}' [LinRL, \S 2.1].  From the nature of the data, one may be able to select a range $\rho$ beyond which one can neglect dependence between two sites (data points $x_i, x_j$ at vertices $v_i, v_j$).  An empirical rule used here is $\rho = \sqrt{8 \nu}/\k$. \\
\i We have stressed the importance of the smoothness parameter $\nu$ (or $\eta$ in {\it 2.2}, {\it 2.5}), yet as we have just seen it is difficult to estimate from data.  This situation is not surprising, and is reminiscent of that in density estimation [Sil].  A discrete distribution can be approximated arbitrarily closely by an arbitrarily smooth one (in any reasonable metric), and vice versa.  So the choice of what smoothness to assume is to be made by the statistician, with the specifics of the particular situation in mind, and the resulting flexibility is an important advantage. \\

\ni {\it 3.2.  Discretization: From Gaussian process to Gaussian Markov random field} \\ 
\i Having chosen the range parameter $\rho$, one can now assume that the data points distant more than $\rho$ apart are {\it conditionally independent given the rest}, that is, by Dempster's theorem, that the corresponding entries in the precision matrix $K$ are zero.  As this will be the case for most entries, the $n \times n$ matrix $K$ is sparse.  This in contrast to the covariance matrix $\Si = K^{-1}$, which will be {\it dense} (most of its entries non-zero).  So matrix operations on $K$ will have a computational cost of $O(n^{3/2})$, while those on $\Si$ have cost $O(n^3)$.  While the `arms race' between increasing computing power and increasing size of data sets is ongoing, such a dramatic difference in computational efficiency is and will remain decisive.  This explains the deliberately dramatic titles of the papers [SimLR1,2] -- which in turn suggested our own. \\
\i In the above, only those edges $e_{ij}$ for which the corresponding $\k_{ij} \neq 0 $ need to be retained (thus the graph will itself be `sparse', reflecting the inter-point distances).  As above, the conditional independence corresponds to the Markov property.  So the continuous Gaussian process has been approximated by a Gaussian (as all distributions are still multivariate normal) Markov random field (GMRF); for background on GMRFs, see Rue and Held [RueH]. \\
\i Although $\Si$ and its inverse $K$ are both symmetric and positive definite, have a Cholesky decomposition etc., they are diametrically opposite computationally as only one is sparse.  Numerical linear algebra for sparse matrices is a subject in its own right; see e.g. Davis et al. [DavRS] and its many references.  We are at liberty to permute the order of the vertices $v_i$, and this can lead to worthwhile improvements in computational efficiency, etc.  Statistical learning with sparsity is also a field in its own right, in which the lasso (least absolute shrinkage and selection operator) plays a central role.  For a monograph treatment, see [HasTW]. \\
\i Discretization of Riemannian manifolds is widely practised; see e.g. Chavel [Cha2, \S 4.4].   There are connections with {\it rough isometries}. \\ 
\i There is much more to say here on the numerics side, particularly on the finite-element method (FEM) and integrated nested Laplace approximation (INLA) [RueMC], but for reasons of space we must refer to the sources cited above, particularly the RSS discussion papers [LinRL], [RueMC].  There is also much to be said about related areas in which sparsity is the key, such as wavelets [Dau] and compressed sensing [FouR].  We mention some of the many application areas.  The first is digitization and reproduction of images: see e.g. [FouR Fig. 1] for a photograph of one of the authors' children, and a near-perfect reconstruction of it using only the 1\% of the wavelet coefficients largest in magnitude.  The second is automatic fingerprint (and more recently, iris) identification systems; see e.g. [MeyC, Fig. 39-42], [RajWMC, Fig. 18-20].  Related applications include automatic number plate recognition for vehicles, bank card recognition for cashless transactions, etc.  \\   

\ni {\bf Dedication and acknowledgement} \\
\i It is a pleasure to dedicate this paper to Mikhailo Iosifovich Yadrenko on his 90th birthday, in acknowledgement of his long and productive scientific career and his many contributions to random fields. \\
\i We thank the editors for their kind invitation to contribute to this special volume, and the referees for their helpful comments. \\

\ni {\bf References} \\
\ni [AbeP] S. Aberg and K. Podg\'orski, A class of non-Gaussian second-order random fields.  {\sl Extremes} {\bf 14} (2011), 187-222. \\
\ni [AndAR] G. E. Andrews, R. Askey and R. Roy, {\sl Special functions}.  Encycl. Math. Appl. {\bf 71}, Cambridge University Press, 1999. \\
\ni [AroS] N. Aronszajn and K. T. Smith, Theory of Bessel potentials I.  {\sl Ann. Inst. Fourier} {\bf 11} (1961), 385-475. \\
\ni [AsaDE] P. Asadi, A. C. Davison and S. Engelke, Extremes on river networks.  {\sl Ann. Appl. Stat.} {\bf 9} (2015), 2023-2050. \\
\ni [AskB] R. Askey and N. H. Bingham, Gaussian processes on compact symmetric spaces.
{\sl Z. Wahrscheinlichkeitstheorie  verw. Geb.} {\bf 37} (1976), 127-143. \\
\ni [BerP] C. Berg and E. Porcu, From Schoenberg coefficients to Schoenberg functions.  {\sl Constructive Approximation} {\bf 45} (2017), 217-241. \\
\ni [Bin1] N. H. Bingham, Random walk on spheres.  {\sl Z. Wahrscheinlichkeitstheorie  verw.
Geb.} {\bf 22} (1972), 169-192. \\
\ni [Bin2] N. H. Bingham, Integral representations for ultraspherical polynomials.  {\sl J.
London Math. Soc.} (2) {\bf 6} (1972), 1-11. \\
\ni [Bin3] N. H. Bingham, Positive definite functions on spheres.  {\sl Proc. Cambridge Phil.
Soc.} {\bf 73} (1973), 145-156. \\
\ni [BinK] N. H. Bingham and  W. J. Krzanowski, The penetration of matrix and linear algebra into multivariate analysis and statistics.  Submitted.\\
%  {\sl Biometrika}.
\ni [BinMS] N. H. Bingham, A. Mijatovi\'c and Tasmin L. Symons, Brownian manifolds, negative type and geo-temporal covariances.  Herbert Heyer Festschrift (ed. D. Applebaum and H. H. Kuo).  {\sl Communications in Stochastic Analysis} {\bf 10} no.4 (2016), 421-432. \\
\ni [BinO] N. H. Bingham and  A. J. Ostaszewski, Extremes and regular variation.  {\sl R. A. Doney Festschrift} (ed. L. Chaumont and G. Peskir);  arXiv:2001.05420. \\
\ni [BinS1] N. H. Bingham and Tasmin L. Symons, Gaussian random fields on sphere and sphere cross line.  {\sl Stochastic Proc. Appl.} (Larry Shepp Memorial Issue); arXiv:1812.02103; {\tt https://doi.org/10.1016/j.spa.2019.08.007}. \\
\ni [BinS2] N. H. Bingham and Tasmin L. Symons, Tasmin L. Symons, Integral representations for ultraspherical polynomials II, arXiv:2101.11809. \\
\ni [BinS3] N. H. Bingham and Tasmin L. Symons, Tasmin L. Symons, Mathematics of Planet Earth, I: Geotremporal covariances.  arXiv:1706.02972. \\
\ni [BinS4] N. H. Bingham and Tasmin L. Symons, Tasmin L. Symons, Mathematics of Planet Earth,
II: The Bochner-Godement theorem for symmetric spaces.  arXiv:1707.05204. \\
\ni [Bog] V. I. Bogachev, {\sl Measure theory}, Volumes 1,2.  Springer, 2007. \\
\ni [BorATMDD] V. Borovitskiy, I. Azangulov, A. Terenin, P. Mostowsky, M. P. Diesenroth and N. Durrande, Mat\'ern Gaussian processes on graphs.  arXiv:2010.15538v3. \\  
\ni [BorTMD] V. Borovitskiy, A. Terenin, P. Mostowsky and M. P. Diesenroth, Mat\'ern Gaussian processes on Riemannian manifolds.  arXiv:2006.10160v2. \\
\ni [Cha1] I. Chavel, {\sl Eigenvalues in Riemannian geometry}.  Academic Press, 1984. \\
\ni [Cha2] I. Chavel, {\sl Riemannian geometry: A modern introduction}.  Cambridge Tracts Math. {\bf 108}, Cambridge University Press, 1993. \\
\ni [Chu] F. R. K. Chung, {\sl Spectral graph theory}.  CBMS Reg. Conf. Ser. Math. {\bf 92}, Amer. Math. Soc., 1997. \\
\ni [CooNN]  D. Cooley, D. Nychka and P. Naveau, Bayesian spatial modelling of extreme precipitation return levels.  {\sl J. Amer. Stat. Assoc.} {\bf 102} (2007), 824-860.  \\
\ni [CraL] H. Cram\'er and M. R. Leadbetter, {\sl Stationary and related stochastic processes.  Sample function properties and their applications}.  Wiley, 1967. \\
\ni [Dav] T. A. Davis, S. Rajamanikam and W. M. Sid-Lakhdar, A survey of direct methods for sparse linear systems.  {\sl Acta Numerica} (2016), 383-566. \\
\ni [Dem1] A. P. Dempster, {\it Elements of continuous multivariate analysis}.  Addison-Wesley, 1969. \\
\ni [Dem2]  A. P. Dempster, Covariance selection.  {\sl Biometrics} {\bf 28} (1972), 157-175. \\ 
\ni [Den] J. Deny, Review of [AroS].  {\sl Mathematical Reviews} {\bf 26}\#1485 (MR0143935).\\
\ni [Dij] G. van Dijk, {\sl Introduction to harmonic analysis and generalized Gelfand pairs}.  De Gruyter Studies in Math. {\bf 36}, Walter de Gruyter, 2009. \\
\ni [Doo] J. L. Doob, {\sl Stochastic processes}.  Wiley, 1953. \\
\ni [Erd] A. Erd\'elyi (ed.), {\sl Higher transcendental functions}, Vol. I-III, McGraw-Hill, 1953. \\
\ni [Fel] W. Feller, {\sl An introduction to probability theory and its applications}, Volume II, end ed.  Wiley, 1971 (1st ed. 1966). \\
\ni [FouR] S. Foucart and H. Rauhut, {\sl A mathematical introduction to compressive sensing}.  Springer, 2013. \\
\ni [Gel] I. M. Gelfand, Spherical functions in symmetric Riemann spaces (Russian).  {\sl Doklady Akac. Nauk SSSR} {\bf 70} (1950), 5-8. \\
\ni [GelV] I. M. Gelfand and N. Ya. Vilenkin, {\sl Generalized functions.  Vol. 4.  Applications of harmonic analysis}.  Academic Press, 1964 (reprinted, AMS Chelsea, 2016). \\
\ni [GenK] M. G. Genton and W. Kleiber, Cross-covariance functions for multivariate geostatistics.  {\sl Stat. Sci.} {\bf 30} (2015), 147-163. \\
\ni [God] R. Godement, Review of [Gel].  {\sl Mathematical Reviews} {\bf 11}, 498b (MR0033832). \\
\ni [GolV] G. H. Golub and C. F. Van Loan, {\sl Matrix computations}, 3rd ed.  The Johns Hopkins University Press, 1996. \\
\ni [GueMP] J. C. Guella, V. A. Menegatto and A. P. Peron, An extension of a theorem of Schoenberg to products of spheres.  {\sl Banach J. Math. Analysis} {\bf 10} (2016), 671-685.\\
\ni [GutG] P. Guttorp and T. Gneiting, Studies in the history of probability and statistics XLIX: On the Mat\'ern correlation family.  {\sl Biometrika} {\bf 93} (2006), 989-995. \\
\ni [HaaRo] L. de Haan and J. de Ronde, Sea and wind: multivariate extremes at work.  {\sl Extremes} {\bf 1} (1998), 7-45. \\
\ni [HasTW] T. Hastie, R. Tibshirani and M. Wainwright, {\sl Statistical learning with sparsity: the lasso and its generalizations}. Monogr. Stat. Appl. Prob. {\bf 143}, Chapman \& Hall/CRC  2015. \\
\ni [Hel1] S. Helgason, {\sl Differential geometry and symmetric spaces}.  Academic Press, 
1962. \\
\ni [Hel2] S. Helgason, {\sl Differential geometry, Lie groups and symmetric spaces}.  Academic Press, 1978. \\
\ni [Hel3] S. Helgason, {\sl Groups and geometric analysis: Integral geometry, invariant differential operators, and spherical functions}.  Academic Press, 1978 (2nd ed., Amer. Math. Soc., 2001). \\
\ni [Hel4] S. Helgason, {\sl Geometric analysis on symmetric spaces}.  Math. Surveys Monog. {\bf 39}, Amer. Math. Soc., 1994 (2nd ed. 2008). \\
\ni [HidKPS] T. Hida, H.-H. Kuo, J. Potthoff and L. Streit, {\sl White noise: An infinite-dimensional calculus}.  Kluwer, 1993. \\
\ni [Hig1] D. Higdon, A process-convolution approach to modelling temperature in the North Atlantic Ocean.  {\sl Environmental and Ecological Statistics} {\bf 5} (1998), 173-190. \\    
\ni [Hig2] D. Higdon, Space and space-time modelling using process convolutions.  {\sl Quantitative methods for current environmental issues} (ed. C. W. Anderson, V. Barnett, P. C. Chatwin and A. H. El-Shaarwi), 37-56.  Springer, 2002. \\
\ni [HorJ] R. A. Horn and C. R. Johnson, {\sl Matrix analysis}.  Cambridge University Press, 1985. \\
\ni [J{\o}r] B. J{\o}rgensen, {\sl Statistical properties of the generalized inverse Gaussian distribution}.  Lecture Notes in Stat. {\bf 9}, Springer, 1982. \\
\ni [Kal1] O. Kallenberg, {\sl Foundations of modern probability}, 2nd ed.  Springer, 2002 (1st ed. 1997). \\
\ni [Kal2] O. Kallenberg, {\sl Random measures, theory and applications}.  Springer, 2017. \\ 
\ni [Lau] S. L. Lauritzen, {\it Graphical models}.  Oxford University Press, 1996.  \\ 
\ni [LinRL] F. Lindgren, H. Rue and T. Lindstr\"om, An explicit link between Gaussian fields and Gaussian Markov random fields.  The SPDE approach (with discussion).  {\sl J. Roy. Stat. Soc. B} {\bf 73} (2011), 423-498. \\
\ni [Mal] A. Malyarenko, {\sl Invariant random fields on spaces with a group action}.  Springer, 2013. \\
\ni [MarP]  D. Marinucci and G. Peccati, {\sl Random fields on the sphere: Representation, limit theorems and cosmological applications}.  London Math. Soc. Lecture Notes Series {\bf 389}, Cambridge University Press, 2011.  \\
\ni [Mat] B. Mat\'ern, {\sl Spatial variation}, 2nd ed.  Springer, 1986 (1st ed. 1960). \\
\ni [MeyC] F. G. Meyer and R. R. Coifman, Brushlets: A tool for directional image analysis and image compression.  {\sl App. Comput. Harm. Anal.} {\bf 4} (1997), 147-187. \\
\ni [PrepS] F. P. Preparata and M. I. Shamos, {\sl Computational geometry: An introduction}.  Springer, 1985. \\
\ni [RajWMC] N. M. Rajpoot, R. G. Wilson, F. G. Meyer and R. R. Coifman, Adaptive wavelet packet basis selection for zerotree image coding.  {\sl IEEE Transations on Image Processing} {\bf 12} (2003), 1460-1471. \\
\ni [Rie] M. Riesz, Int\'egrales de Riemann-Liouville et potentiels.  {\sl Acta Szeged Sect. Math.} {\bf 9} (1937-38), 1-42; reprinted, Paper 42 p.482-523, in {\sl Collected papers} (ed. L. Garding and L. H\"ormander).  Springer, 1988. \\
\ni [RodD] A. Rodrigues and P. J. Diggle, A class of convolution-based models for spatio-temporal processes with non-separable covariance structure.  {\sl Scand. J. Stat.} {\bf 37} (2010), 553-567. \\
\ni [Rud] W. Rudin, {\sl Functional analysis}, 2nd ed.  McGraw-Hill, 1991 (1st ed. 1973). \\
\ni [RueH] H. Rue and L. Held, {\sl  Markov random fields: Theory and applications.}  Monographs in Statistics and Applied Probability {\bf 104}, Chapman \& Hall/CRC, 2005. \\
\ni [RueMC] H. Rue, S. Martino and N. Chopin, Approximate Bayesian inference for latent Gaussian models by using integrated nested Laplace approximation (with discussion).  {\sl J. Roy. Stat. Soc. B} {\bf 71} (2009), 319-392. \\
\ni [Sch] I. Schur, Bemerkungen zur Theorie der beschr\"ankten Bilinearformen mit unendlich vielen Ver\"anderlichen.  {\sl J. Reine Angew. Math.} {\bf 140} (1911), 1-28. \\
\ni [ShaW] P. Sharkey and H. C. Winter, A Bayesian spatial hierarchical model for extreme precipitation in Great Britain.  {\sl Environmetrics} {\bf 30} (2019), e2529.  \\
\ni [Sil] B. W. Silverman, {\sl Density estimation for statistics and data analysis}.  Chapman \& Hall, 1986. \\
\ni [SimLR1] D. Simpson, F. Lindgren and H. Rue, In order to make spatial statistics computationally feasible, we need to forget about the covariance function.  {\sl Environmetrics} {\bf 23}(1) (Special Issue on Spatio-Temporal Stochastic Modelling) (2012), 65-74. \\
\ni [SimLR2] D. Simpson, F. Lindgren and H. Rue, Beyond the valley of the covariance function.  {\sl Stat. Sci.} {\bf 30} (2015), 164-166. \\
\ni [Ste] E. M. Stein, {\sl Singular integrals and differentiability properties of functions}.  Princeton University Press, 1970. \\
\ni [SteW] E. M. Stein and G. Weiss, {\sl Introduction to Fourier analysis on Euclidean spaces}.  Princeton University Press, 1971. \\
\ni [Stei] M. L. Stein, {\sl Interpolation of spatial data: Some theory for kriging}.  Springer, 1999. \\
\ni [Sze] G. Szeg\H{o}, {\sl Orthogonal polynomials}.  AMS Colloquium Publications {\bf XXIII}, Amer. Math. Soc., 1959. \\
\ni [Tay] M. E. Taylor, {\sl Pseudodifferential operators}.  Princeton University Press, 1981. \\
\ni [Wal] J. B. Walsh, An introduction to stochastic partial differential equations.  {\sl Ecole d'Et\'e de Saint-Flour XIV}, 265-439; {\sl Lecture Notes in Math.} {\bf 1180}, Springer, 1984. \\ 
\ni [Wat] G. N. Watson, {\sl A treatise on the theory of Bessel functions}, 2nd ed.  Cambridge University Press. 1944 (1st ed. 1922). \\
\ni [Whi1] P. Whittle, On stationary processes in the plane.  {\sl Biometrika} {\bf 41} (1954), 434-449. \\
\ni [Whi2] P. Whittle, Stochastic processes in several dimensions.  {\sl Bull. Internat. Stat. Inst.} {\bf 40} (1963), 974-994. \\
\ni [Whi3] P. Whittle, {\sl Systems in stochastic equilibrium}.  Wiley, 1986. \\
\ni [Wil1] D. Williams, {\sl Probability with martingales}.  Cambridge University Press, 1991. \\
\ni [Wil2] D. Williams, {\sl Weighing the odds: A course in probability and statistics}.  Cambridge University Press, 2001. \\
\ni [Wol1] Wolf, J. A.: {\sl Spaces of constant curvature}.  Amer. Math. Soc., 1967 (6th ed. 2011). \\
\ni [Wol2] Wolf, J. A.: {\sl Harmonic analysis on commutative spaces}.  Amer. Math. Soc., 2007.\\ 
\ni [Yad] M. I. Yadrenko, {\sl Spectral theory of random fields}.  Optimization Software, 1983.\\

\ni Mathematics Department, Imperial College, London SW7 2AZ, UK; n.bingham@ic.ac.uk \\

\ni Telethon Kids Institute, 15 Hospital Avenue, Perth, WA 6009, Australia; \\ tasmin.symons@telethonkids.org.au    

\end{document}